\documentclass[runningheads]{llncs}

\usepackage{epsfig}
\usepackage{url}
\usepackage{amsmath}
\usepackage{amssymb}

\newcommand{\NP}{\mathcal{NP}}
\newcommand{\order}[1]{\mathcal{O}\left({#1}\right)}

\newcommand{\orderbig}[1]{\mathcal{O}\big({#1}\big)}
\newcommand{\coNP}{\text{co}\mathcal{NP}}

\newcommand{\nodes}[1]{V_{#1}}
\newcommand{\edges}[1]{E_{#1}}
\newcommand{\graph}[1]{G_{#1}}
\newcommand{\lattice}[1]{L_{#1}}
\newcommand{\Hsum}[1]{\mathcal{H}\left({#1}\right)}
\newcommand{\Htwosum}[1]{\mathcal{H}_2\left({#1}\right)}

\newcommand{\auxG}[1]{G^{\star}_{#1}}
\newcommand{\twofac}[1]{\mathcal{C}_{#1}}

\newcommand{\neighbors}[1]{\Gamma\left({#1}\right)}
\newcommand{\spannedFace}[1]{F\left({#1}\right)}
\newcommand{\dual}[1]{{#1}^{\star}}
\newcommand{\boundary}[1]{\partial{#1}}

\begin{document}

\pagestyle{headings}

\mainmatter

\title{Reconstructing a Simple Polytope\\%
       from its Graph}

\titlerunning{Reconstructing a Simple Polytope from its Graph}

\author{Volker Kaibel%
  \thanks{Supported by the Deutsche
    Forschungsgemeinschaft, FOR~413/1--1 (Zi~475/3--1).
    }
}

\authorrunning{Volker Kaibel}

\institute{TU~Berlin, MA~6--2\\
Stra\ss e des 17. Juni~136\\
10623~Berlin, Germany\\
\email{kaibel@math.tu-berlin.de}\\
\url{http://www.math.tu-berlin.de/~kaibel}}

\maketitle

\begin{abstract}
  Blind and Mani~\cite{BM87} proved that the entire combinatorial
  structure (the vertex-facet incidences) of a simple convex polytope
  is determined by its abstract graph. Their proof is not
  constructive. Kalai~\cite{Kal88} found a short, elegant, and
  algorithmic proof of that result. However, his algorithm has always
  exponential running time. We show that the problem to reconstruct
  the vertex-facet incidences of a simple polytope~$P$ from its graph
  can be formulated as a combinatorial optimization problem that is
  strongly dual to the problem of finding an abstract objective
  function on~$P$ (i.e., a shelling order of the facets of the dual
  polytope of~$P$). Thereby, we derive polynomial certificates for
  both the vertex-facet incidences as well as for the abstract
  objective functions in terms of the graph of~$P$. The paper is a
  variation on joint work with Michael Joswig and Friederike
  K\"orner~\cite{JKK01}.
\end{abstract}

%                                                         
% Introduction                                            
%                                                         
\section{Introduction}
\label{sec:intro}

The \emph{face lattice}~$\lattice{P}$ of a (convex) polytope~$P$ is
any lattice that is isomorphic to the lattice formed by the set of all
faces of~$P$ (including~$\varnothing$ and~$P$ itself), ordered by
inclusion. It is well-known to be determined by the \emph{vertex-facet
  incidences} of~$P$, i.e., by any graph that is isomorphic to the
bipartite graph whose nodes are the vertices and the facets of~$P$,
where the edges are defined by the pairs $\{v,f\}$ of vertices~$v$ and
facets~$f$ with $v\in f$. In lattice theoretic terms, $\lattice{P}$ is
a ranked, atomic, and coatomic lattice, and thus, the sub-poset
formed by its atoms and coatoms already determines the whole lattice.
Actually, one can compute~$\lattice{P}$ from the vertex-facet
incidences of~$P$ in $\order{\eta\cdot\alpha\cdot\lambda}$ time,
where~$\eta$ is the minimum of the number of vertices and the number
of facets, $\alpha$ is the number of vertex-facet incidences,
and~$\lambda$ is the total number of faces of~$P$~\cite{KP01}.

\begin{figure}[ht]
  \begin{center}
    \epsfig{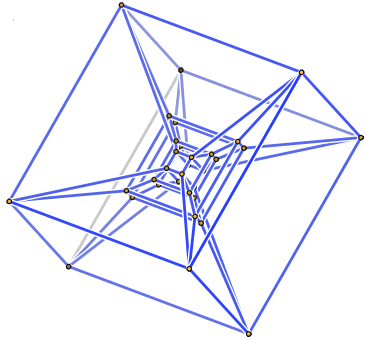}
    \caption{A Schlegel-diagram (projection onto one facet) of a four-dimensional polytope with the graph of a five-dimensional cube, found by Joswig \&\ Ziegler~\cite{JZ00}.}
    \label{fig:C45}
  \end{center}
\end{figure}

The \emph{graph}~$\graph{P}=(\nodes{P},\edges{P})$ of a polytope~$P$
is any graph that is isomorphic to the graph whose nodes are the
vertices of~$P$, where two nodes are adjacent if and only if the
convex hull of the corresponding two vertices is a one-dimensional
face of~$P$. Phrased differently, $\graph{P}$ is the graph defined on
the rank one elements of~$\lattice{P}$, where two rank one elements
are adjacent if and only if they are below a common rank two element.
While the vertex-facet incidences completely determine the
face-lattice of any polytope, the graph of a polytope in general does
not encode the entire combinatorial structure. This can be seen, e.g.,
from the examples of the cut polytope associated with the complete
graph on~$n$ nodes and the $\binom{n}{2}$-dimensional cyclic polytope
with $2^{n-1}$ vertices, which both have complete graphs. Another
example is the four-dimensional polytope shown in Fig.~\ref{fig:C45}
whose graph is isomorphic to the graph of the five-dimensional cube.

Actually, in all dimensions less than four such ambiguities cannot
occur. For one- or two dimensional polytopes this is obvious, and for
three-dimensional polytopes it follows from Whitney's
theorem~\cite{Whi32} saying that every $3$-connected planar graph has
a unique (up to reflection) plane embedding.

A $d$-dimensional polytope~$P$ is \emph{simple} if every vertex of~$P$
is contained in precisely~$d$ facets, which is equivalent
to~$\graph{P}$ being $d$-regular (the polytope is non-degenerate in
terms of Linear Programming).  Every face of a simple polytope is
simple as well.  None of the examples showing that it is, in general,
impossible to reconstruct the face lattice of a polytope from its
graph, is simple.

In fact, Blind and Mani~\cite{BM87} proved in 1987 that the face
lattice of a simple polytope is determined by its graph.  Their
proof (which we sketch in Sect.~\ref{sec:BM}) is not constructive and
crucially relies on the topological concept of homology. In 1988,
Kalai~\cite{Kal88} found a short and elegant proof (reviewed in
Sect.~\ref{sec:K}) that does only use elementary geometric and
combinatorial reasoning with the main advantage of being algorithmic.
However, the running time of the method that can be devised from it is
exponential in the size of the graph.

Perles conjectured in the 1970's (see~\cite{Kal88}) that for a
$d$-dimensional simple polytope~$P$ every subset $F\subset\nodes{P}$
that induces a $(d-1)$-regular, connected, and non-separating subgraph
of~$\graph{P}$ corresponds to the vertex set of a facet of~$P$. A
proof of this conjecture would have lead immediately to a polynomial
time algorithm that, given the graph~$\graph{P}=(\nodes{P},\edges{P})$
of a simple polytope~$P$, decides for a set of subsets of~$\nodes{P}$
if it corresponds to the set of vertex sets of facets of~$P$.
However, Haase and Ziegler~\cite{HZ01} recently disproved Perles'
conjecture. They found a four-dimensional simple polytope whose graph
has a $3$-regular, non-separating, and even $3$-connected induced
subgraph that does not correspond to any facet.

Refining ideas from Kalai's proof (Sect.~\ref{sec:WO}), we show that
the problem of reconstructing the vertex-facet incidences of a simple
polytope~$P$ from its graph~$\graph{P}$ can be formulated as a
combinatorial optimization problem that has a well-stated strongly
dual problem (Sect.~\ref{sec:GC}). The optimal solutions to this dual
problem are certain orientations of~$\graph{P}$ (induced by ``abstract
objective functions'') that are important also in different contexts.
In particular, we provide short certificates for both the vertex-facet
incidences of a simple polytope and for the abstract objective
functions in terms of~$\graph{P}$. We conclude with some remarks on
the complexity status of the problem to decide whether a claimed
solution to the reconstruction problem indeed are the vertex-facet
incidences of the respective polytope in Sect.~\ref{sec:rem}.

The material presented here has evolved from joint work with Michael
Joswig and Friederike K\"orner~\cite{JKK01}.  The basic ideas and
results are the same in both papers. However, the concept of a
``facoidal system of walks'' is newly introduced here. It differs from
the corresponding notion of a ``$2$-system'' (introduced
in~\cite{JKK01}) with the main effect that one knows how to compute
efficiently some facoidal system of walks from the graph ~$\graph{P}$
of a simple polytope~$P$ (see Proposition~\ref{prop:facoidal}), while it is
unclear how to find a $2$-system from~$\graph{P}$ efficiently.
Furthermore, the proof of Theorem~\ref{thm:2facesAOF} we give here is
different from the corresponding proof in~\cite{JKK01}. Finally, the
complexity theoretic statement of Corollary~\ref{cor:steinitz} does
not appear in~\cite{JKK01}.

For all notions from the theory of polytopes that we use without
(sufficient) explanations, we refer to Ziegler's book~\cite{Zie95}.  We
use the terms \emph{$d$-polytope} and \emph{$k$-face} for
$d$-dimensional polytopes and $k$-dimensional faces, respectively.
Often we will identify a face~$F$ of a simple polytope~$P$ with the
subset of nodes of~$\graph{P}$ that corresponds to the vertex set
of~$F$.  Whenever we talk about ``polynomial time'' or ``efficient''
this refers to the size of the graph~$\graph{P}$ of the respective
(simple) polytope~$P$.

%                                                         
% The Theorem of Blind and Mani                             
%                                                         
\section{The Theorem of Blind and Mani}
\label{sec:BM}

Blind and Mani~\cite{BM87} proved their theorem in the dual setting,
i.e., for simplicial rather than for simple polytopes. Nevertheless,
we sketch parts of their proof in terms of simple polytopes here. The
starting point is the observation that, while a priori it is by no
means clear if the graph~$\graph{P}=(\nodes{P},\edges{P})$ of a simple
polytope~$P$ determines the face lattice of~$P$, it is easy to see
that the $2$-faces of~$P$ (as subsets of~$\nodes{P}$) carry the entire
information on the combinatorial structure of~$P$.

Let~$P$ be a simple polytope. For a node $v\in\nodes{P}$
of~$\graph{P}$ denote by $\neighbors{v}\subset\nodes{P}$ the subset of
nodes that are adjacent to~$v$ (the \emph{neighbors} of~$v$).  For any
$k$-element subset $S\subset\neighbors{v}$ there is a $k$-face of~$P$
that contains the vertices that correspond to $S\cup\{v\}$ (and no
vertex that corresponds to a node in $\neighbors{v}\setminus S$).  We
call the subset $\spannedFace{S\cup\{v\}}\subseteq\nodes{P}$ of nodes
corresponding to the vertices of that face the \emph{$k$-face spanned
  by $S\cup\{v\}$}.

For an edge~$e=\{v,w\}\in\edges{P}$ let
$\Psi_{(v,w)}:\neighbors{v}\setminus\{w\}\longrightarrow\neighbors{w}\setminus\{v\}$
be the map that assigns to each subset
$S\subset\neighbors{v}\setminus\{w\}$ the subset
$T\subset\neighbors{w}\setminus\{v\}$ with
$\spannedFace{S\cup\{v,w\}}=\spannedFace{T\cup\{w,v\}}$.  The
maps~$\Psi_{(v,w)}$ are cardinality preserving bijections, where
$\Psi_{(w,v)}$ is the inverse of~$\Psi_{(v,w)}$.

\begin{proposition}
  \label{prop:2facesvfi}
  Let~$P$ be a simple polytope.  For each $e=\{v,w\}\in\edges{P}$
  and $S\subseteq\neighbors{v}\setminus\{w\}$ we have
  $$
  \Psi_{(v,w)}(S)=\overline{\Psi_{(v,w)}(\overline{S})}
  $$
  (where $\overline{U}$ is the respective complement of the set~$U$).
\end{proposition}

\begin{proof}
  This follows from the fact that we have
  $\Psi_{(v,w)}(S_1\cap
  S_2)=\Psi_{(v,w)}(S_1)\cap\Psi_{(v,w)}(S_2)$ for all
  $S_1,S_2\subseteq\neighbors{v}\setminus\{w\}$.
\end{proof}

With the notations of Proposition~\ref{prop:2facesvfi} denote by
$\Psi^k_{(v,w)}$ the restriction of the map~$\Psi^k_{(v,w)}$ to the
$(k-1)$-element subsets of~$\neighbors{v}\setminus\{w\}$. There are
quite obvious algorithms that compute from the maps~$\Psi^k_{(v,w)}$,
$\{v,w\}\in\edges{P}$, the $k$-faces of~$P$ (as subsets
of~$\nodes{P}$), and vice versa, in polynomial time in the
number~$f_k(P)$ of~$k$-faces of the simple $d$-polytope~$P$. Since
both~$f_2(P)$ as well as $f_{d-1}(P)$ are bounded polynomially in the
size of~$\graph{P}$, the following result follows.

\begin{corollary}
  \label{cor:2facesvfi}
  There are polynomial time algorithms that, given the
  graph~$\graph{P}$ of a simple $d$-polytope~$P$, compute the set of
  facets of~$P$ from the set of $2$-faces of~$P$ (both viewed as sets
  of subsets of~$\nodes{P}$), and vice versa.
\end{corollary}

For the rest of this section let~$P_1$ and~$P_2$ be two simple
polytopes, and let~$g:\nodes{P_1}\longrightarrow\nodes{P_2}$ be an
isomorphism of the graphs~$\graph{P_1}=(\nodes{P_1},\edges{P_1})$
and~$\graph{P_2}=(\nodes{P_2},\edges{P_2})$ of~$P_1$ and~$P_2$,
respectively (i.e., $g$ is an in both directions edge preserving
bijection). 

The core of Blind and Mani's paper~\cite{BM87} is the following
result.

\begin{proposition}
\label{prop:2face22face}
  The graph isomorphism~$g$ maps every cycle in~$\graph{P_1}$ that
  corresponds to a $2$-face of~$P_1$ to a cycle in~$\graph{P_2}$ that
  corresponds to a $2$-face of~$P_2$.
\end{proposition}

Blind and Mani's proof proceeds in the dual setting, i.e., in terms of
the boundary complexes~$\boundary{\dual{P_1}}$
and~$\boundary{\dual{P_2}}$ of the simplicial dual
polytopes~$\dual{P_1}$ and~$\dual{P_2}$ of~$P_1$ and~$P_2$,
respectively. The strategy is to show that, if some cycle
in~$\graph{P_1}$ corresponding to a $2$-face of~$P_1$ was mapped to
some cycle in~$\graph{P_2}$ that does not correspond to any $2$-face
of~$P_2$, then a certain sub-complex of~$\boundary{\dual{P_2}}$ would
have a certain non-vanishing (reduced) homology group. They complete
their proof of Proposition~\ref{prop:2face22face} by showing that the
respective homology group, however, is zero. The key ingredient they
use to prove this is the following. For each face~$F$ of~$\dual{P_1}$
there is a \emph{shelling order} of the facets of~$\dual{P_1}$ (i.e.,
an ordering satisfying certain convenient topological properties,
which, however, can be expressed completely combinatorially, see
Sect.~\ref{sec:K}) in which the facets containing~$F$ come first.

From Proposition~\ref{prop:2face22face} one can deduce that the graph
isomorphism~$g$ actually induces a bijection between the cycles
in~$\graph{P_1}$ that correspond to $2$-faces of~$P_1$ and the cycles
in~$\graph{P_2}$ that correspond to $2$-faces of~$P_2$. Once this is
established, Proposition~\ref{prop:2facesvfi} yields the following
result.

\begin{theorem}[Blind \&\ Mani~\cite{BM87}]
  \label{thm:BM}
  Every isomorphism between the graphs~$\graph{P_1}$ and~$\graph{P_2}$
  of two simple polytopes~$P_1$ and~$P_2$, respectively, induces an
  isomorphism between the vertex-facet incidences of~$P_1$
  and~$P_2$. In particular, the graph of a simple polytope determines
  its entire face lattice.
\end{theorem}

%                                                         
% Kalai's Constructive Proof                              
%                                                         
\section{Kalai's Constructive Proof}
\label{sec:K}

Kalai realized that the existence of shelling orders as exploited by
Blind and Mani can be used directly in order to devise a simple proof
which does not rely on any topological notions like
homology~\cite{Kal88}.  He formulated his proof in the original
setting, i.e., for simple polytopes, where the notion corresponding to
``shelling'' is called ``abstract objective function.''

From now on, let~$P$ be a simple $d$-polytope with~$n$ vertices. For
simplicity of notation, we will identify each face of~$P$ not only
with the corresponding subset of~$\nodes{P}$, but also with the
corresponding induced subgraph of~$\graph{P}$. Furthermore, by saying
that $w\in W\subset\nodes{P}$ is a sink of~$W$ we mean that~$w$ is a
sink of the orientation induced on the subgraph of~$\graph{P}$ that is
induced by~$W$.

\begin{definition}
  Every bijection $\varphi:\nodes{P}\longrightarrow\{1,\dots,n\}$
  induces an acyclic orientation~$\mathcal{O}_{\varphi}$ of the
  graph~$\graph{P}$ of~$P$, where an edge is directed from its larger
  end-node to its smaller end-node (with respect to~$\varphi$).  The
  map~$\varphi$ is called an \emph{abstract objective function}
  (\emph{AOF}) if~$\mathcal{O}_{\varphi}$ has a unique sink in every
  non-empty face of~$P$ (including~$P$ itself). Such an orientation
  of~$\graph{P}$ is called an \emph{AOF-orientation}.
\end{definition}

The inverse orientation of an AOF-orientation is an AOF-orientation as
well (this follows, e.g., from Theorem~\ref{thm:2facesAOF}). Thus,
every AOF-orientation also has a unique source in every non-empty
face.

From the fact that the simplex algorithm works correctly (on every
face) one easily derives that every linear function that assigns
pairwise different values to the vertices of~$P$ induces an
AOF-orientation (this is a consequence of the convexity of the faces).
From this observation, the following fact follows (which is dual to
the existence of the shelling orders required in Blind and Mani's
proof).

\begin{lemma}
  \label{lem:manyAOF}
  Let~$W\subset\nodes{P}$ be any face of~$P$. There is an
  AOF-orientation of~$\graph{P}$ for which~$W$ is \emph{terminal},
  i.e., no edge in the cut defined by~$W$ is directed from~$W$
  to~$\nodes{P}\setminus W$.
\end{lemma}

In a sense, this statement can be reversed.

\begin{lemma}
  \label{lem:AOFterm}
  Let~$W\subset\nodes{P}$ be a set of nodes inducing a~$k$-regular
  connected subgraph of~$\graph{P}$, and let~$\mathcal{O}$ be an
  AOF-orientation for which~$W$ is terminal.  Then~$W$ is a~$k$-face
  of~$P$.
\end{lemma}

\begin{proof}
  Since~$\mathcal{O}$ is acyclic, it has a source~$s$ in~$W$.
  Let~$w_1,\dots,w_k\in W$ be the neighbors of~$s$ in~$W$, and
  let~$F:=\spannedFace{\{t,w_1,\dots,w_k\}}\subset\nodes{P}$ be the
  $k$-face of~$P$ that is spanned by~$t,w_1,\dots,w_k$.
  Since~$\mathcal{O}$ has unique sources on non-empty faces, $s\in
  W\cap F$ must be the unique source of~$F$. By the acyclicity
  of~$\mathcal{O}$ there hence is a monotone path from~$s$ to every
  node in~$F$. Since~$W$ is terminal this implies $F\subseteq W$.
  Because both~$F$ and~$W$ induce $k$-regular connected subgraphs
  of~$\graph{P}$, $F=W$ follows.
\end{proof}

Lemma~\ref{lem:manyAOF} and Lemma~\ref{lem:AOFterm} imply that one can
compute the vertex-facet incidences of~$P$, provided that one knows
all AOF-orientations of~$\graph{P}$. Kalai's crucial discovery is that
one can compute the AOF-orientations just from~$\graph{P}$ (i.e.,
without explicitly knowing the faces of~$P$).

\begin{definition}
  For an orientation~$\mathcal{O}$ of~$\graph{P}$
  let~$h_k(\mathcal{O})$ be the number of nodes with in-degree~$k$. The
  number
  $$
  \Hsum{\mathcal{O}}:=\sum_{k=0}^{d}h_k(\mathcal{O})\cdot 2^k
  $$
  is called the \emph{$\mathcal{H}$-sum} of~$\mathcal{O}$.
\end{definition}

Since every subset of neighbors of a vertex~$v$ of~$P$ together
with~$v$ span a face of~$P$ containing no other neighbors of~$v$, one
finds (by double-counting) that
\begin{equation}
  \label{eq:Hsum}
  \Hsum{\mathcal{O}}=\sum_{F\text{ face of }P}\big(\#\text{ sinks
    of~$\mathcal{O}$ in~$F$\big)}   
\end{equation}
is the total number of sinks induced by~$\mathcal{O}$ on faces of~$P$.
Consequently, since every acyclic orientation has at least one sink in
every non-empty face, we have the following characterization.

\begin{lemma}
  \label{lem:minHsum}
  An orientation~$\mathcal{O}$ of~$\graph{P}$ is an AOF-orientation if
  and only if it is acyclic and has minimal $\mathcal{H}$-sum among
  all acyclic orientations of~$\graph{P}$ (which then equals the
  number of non-empty faces of~$P$).
\end{lemma}

Thus, by enumerating all $2^{\frac{d\cdot n}{2}}=\sqrt{2}^{d\cdot n}$
orientations of~$\graph{P}$ one can find all AOF-orientations
of~$\graph{P}$.

\begin{theorem}[Kalai~\cite{Kal88}]
  \label{thm:kalai}
  There is an algorithm that computes the vertex-facet incidences of a
  simple $d$-polytope with~$n$ vertices from its graph in
  $\orderbig{\sqrt{2}^{d\cdot n}}$ steps.
\end{theorem}

%                                                         
% Walks and Orientations                                  
%                                                         
\section{Walks and Orientations}
\label{sec:WO}

In this section, we refine the ideas of Kalai's proof and combine them
with the observation (exploited by Blind and Mani) that it suffices to
identify the $2$-faces from the graph of a simple polytope, even with
respect to the question for polynomial time reconstruction algorithms
(see Corollary~\ref{cor:2facesvfi}). Let us start with a result that
emphasizes the importance of the $2$-faces even more. The result was
known for cubes~\cite{HSLdW88}; for three-dimensional simple polytopes
it was independently proved by Develin~\cite{Dev00}). For general
simple polytopes it seems that it was assumed to be false
(see~\cite[Ex.~8.12~(iv)]{Zie95}).

\begin{theorem}
\label{thm:2facesAOF}
  An acyclic orientation~$\mathcal{O}$ of the graph~$\graph{P}$ of a
  simple polytope~$P$ is an AOF-orientation if and only if it has a
  unique sink on every $2$-face of~$P$.
\end{theorem}

\begin{proof}
  The ``only if'' part is clear by definition. For the ``if'' part,
  let $\varphi:\nodes{P}\longrightarrow\{1,\dots,n\}$ be a bijection
  inducing an acyclic orientation~$\mathcal{O}=\mathcal{O}_{\varphi}$
  that has a unique sink in every $2$-face of~$P$. Suppose there is a
  face~$F$ of~$P$ in which~$\mathcal{O}$ has two sinks~$t_1,t_2\in F$
  ($t_1\not=t_2$). We might assume~$F=P$ (because~$F$ itself is a
  simple polytope with every $2$-face of~$F$ being a $2$-face of~$P$).
  
  Since~$\graph{P}$ is connected, there is a path in~$\graph{P}$
  connecting~$t_1$ and~$t_2$.  Let~$\Pi\not=\varnothing$ be
  the set of all these paths. For every $\pi\in\Pi$ we denote
  by~$\mu(\pi)$ the maximal $\varphi$-value of any node in~$\pi$.
  Let~$\pi_{\min}\in\Pi$ be a path with minimal $\mu$-value
  among all paths in~$\Pi$, where~$v\in\nodes{P}$ is the node
  in~$\pi_{\min}$ with $\varphi(v)=\mu(\pi_{\min})$ (see
  Fig.~\ref{fig:2facesAOF}).
  \begin{figure}[ht]
    \begin{center}
      \epsfig{figure=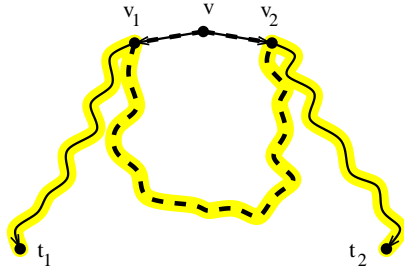,height=3.5cm}
      \caption{Illustration of the proof of
        Theorem~\ref{thm:2facesAOF}. The fat grey path is the one
        yielding the contradiction.}
      \label{fig:2facesAOF}
    \end{center}
  \end{figure}
  
  Obviously, $v$ is a source in the path~$\pi_{\min}$ (in particular,
  $v\not\in\{t_1,t_2\}$). Let~$C$ be the $2$-face spanned by~$v$ and
  its two neighbors~$v_1$ and~$v_2$ in~$\pi_{\min}$. Since~$v$ is the
  unique source of~$\mathcal{O}$ in~$C$, $v$ has the largest
  $\varphi$-value among all nodes in the union~$U$ of~$C$
  and~$\pi_{\min}$. But $U\setminus\{v\}$ induces a connected subgraph
  of~$\graph{P}$ containing both~$t_1$ and~$t_2$, which contradicts
  the minimality of $\mu(\pi_{min})$.
\end{proof}

From now on let, again, $P$ be a simple polytope.  The ultimate goal
is to find the system of cycles in the graph~$\graph{P}$ that
corresponds to the set of $2$-faces of~$P$. However, we even do not
know how to \emph{prove} or \emph{disprove} efficiently that a given
system of cycles actually is the one we are searching for.  We now
define more general systems having the property that one can at least
generate one of them in polynomial time (which in general, of course,
will not be the desired one), and among which the one corresponding to
the set of $2$-faces of~$P$ can be characterized using
AOF-orientations.

\begin{definition}
  \begin{enumerate}
  \item A sequence $W=(w_0,\dots,w_{l-1})$ (with $l\geq 3$) of nodes
    in~$\graph{P}$ is called a \emph{closed smooth walk} in~$\graph{P}$, if
    $\{w_i,w_{i+1}\}$ is an edge of~$\graph{P}$ and
    $w_{i-1}\not=w_{i+1}$ for all $i$ (where, as in the following, all
    indices are taken modulo~$l$).  Note that the~$w_i$ need not be
    pairwise disjoint. We will identify two closed smooth walks if they
    differ only by a cyclic shift and/or a ``reflection'' of their
    node sequences.
  \item A set~$\mathcal{W}$ of closed smooth walks in~$\graph{P}$ is a
    \emph{facoidal system of walks} if for every
    triple~$v,v_1,v_2\in\nodes{P}$ ($v_1\not= v_2$) such that
    both~$v_1$ and~$v_2$ are neighbors of~$v$ there is a \emph{unique}
    closed smooth walk $(w_0,\dots,w_{l-1})\in\mathcal{W}$ with
    $(w_{i-1},w_i,w_{i+1})\in\{(v_1,v,v_2),(v_2,v,v_1)\}$ for
    some~$i$, which is also required to be unique.
  \end{enumerate}
\end{definition}

\begin{figure}[ht]
  \begin{center}
    \epsfig{figure=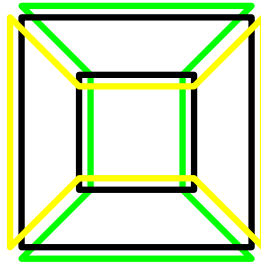,height=3.5cm}
    \caption{A facoidal system of four walks in the graph of the three-dimensional cube.}
    \label{fig:otherFacoidal}
  \end{center}
\end{figure}

The system of $2$-faces of~$P$ yields a uniquely determined (recall
the identifications mentioned in part~1 of the definition) facoidal
system of walks in~$\graph{P}$, which is denoted by~$\twofac{P}$.  In
general, there are many other facoidal systems of walks (see
Fig.~\ref{fig:otherFacoidal}).

%The \emph{line graph} $\lineG{G}$ of a graph~$G=(V,E)$ is the graph
%whose nodes are the edges~$E$ of~$G$, where two (different) edges
%$e_1,e_2\in E$ of~$G$ form an edge~$\lambda$ of~$\lineG{G}$ if they
%have a common end-node~$v(\lambda)$ (i.e., $e_1\cap
%e_2=\{v(\lambda)\}$). Let $\auxG{P}$ be the subgraph
%of~$\lineG{\lineG{\graph{P}}}$, which contains only those
%edges~$\{\lambda_1,\lambda_2\}$ of~$\lineG{\lineG{\graph{P}}}$ with
%$v(\lambda_1)\not=v(\lambda_2)$ (see Fig.~\ref{fig:llgraph}).

\begin{figure}[ht]
  \begin{center}
    \epsfig{figure=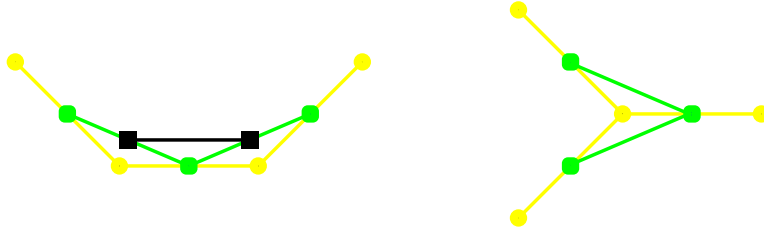,height=3cm}
    \caption{The left constellation gives rise to an edge of~$\auxG{P}$, while the right one does not.}
    \label{fig:llgraph}
  \end{center}
\end{figure}

For each path~$\lambda$ in~$\graph{P}$ of length two denote
by~$v(\lambda)$ the inner node of~$\lambda$. Let~$\auxG{P}$ be the
graph defined on the paths of length two in~$\graph{P}$, where two
paths~$\lambda_1$ and~$\lambda_2$ are adjacent if and only if they
share a common edge and $v(\lambda_1)\not=v(\lambda_2)$ holds (see
Fig.~\ref{fig:llgraph}).

A \emph{$2$-factor} in a graph~$G$ is a set of (not self-intersecting)
cycles in~$G$ such that every node is contained in a unique cycle.
Checking whether a graph has a $2$-factor and finding one (if it
exists) can be reduced (by a procedure due to Tutte~\cite{Tut54}) to
searching for a perfect matching in a related graph (which can be
performed in polynomial time by Edmonds' algorithm~\cite{Edm65a}).

\begin{proposition}
\label{prop:facoidal}
For simple polytopes~$P$,
  \begin{enumerate}
  \item there is a (polynomial time computable) bijection between the
    facoidal systems of walks in~$\graph{P}$ and the $2$-factors of
    $\auxG{P}$,
  \item checking whether a given set of node-sequences in~$\graph{P}$
    is a facoidal system of walks can be done in polynomial time, and
  \item one can find a facoidal system of walks in~$\graph{P}$ in
    polynomial time.
  \end{enumerate}
\end{proposition}

\begin{proof}
  Part~2 is obvious, part~3 follows from part~1 by Tutte's
  reduction~\cite{Tut54} and Edmonds algorithm~\cite{Edm65a}, and
  part~1 is readily obtained from the definitions.
\end{proof}

Proposition~\ref{prop:facoidal} shows that facoidal systems of walks have
quite convenient algorithmic properties. However, they become useful
only due to the fact that the system~$\twofac{P}$ corresponding to the
$2$-faces of~$P$ can be well-characterized among them, as we will
demonstrate next.

\begin{definition}
  Let~$\mathcal{O}$ be any orientation of~$\graph{P}$. 
  \begin{enumerate}
  \item The \emph{$\mathcal{H}_2$-sum} of~$\mathcal{O}$ is defined as
    $$
    \Htwosum{\mathcal{O}}:=\sum_{k=0}^{d}h_k(\mathcal{O})\cdot \binom{k}{2}\enspace.
    $$
  \item A closed smooth walk~$(w_0,\dots,w_{l-1})$ in~$\graph{P}$ has
    a \emph{sink} (\emph{source}, respectively) at position~$i$ (with
    respect to the orientation~$\mathcal{O}$), if the
    edges~$\{w_i,w_{i-1}\}$ and~$\{w_i,w_{i+1}\}$ both are directed
    towards (away from, respectively)~$w_i$.
  \end{enumerate}
\end{definition}

The following follows immediately from the definitions.

\begin{lemma}
  \label{lem:H2sum}
  For every orientation~$\mathcal{O}$ of~$\graph{P}$ the
  sum~$\Htwosum{\mathcal{O}}$ equals the total number of sinks (with
  respect to~$\mathcal{O}$) in \emph{every} facoidal system of walks
  in~$\graph{P}$.
\end{lemma}

Now we can formulate and prove the main result of this section (where
$f_2(P)$ denotes the number of $2$-faces of~$P$).

\begin{theorem}
  \label{thm:walkOrient}
  Let~$P$ be a simple polytope, $\mathcal{W}$ a facoidal system of
  walks in~$\graph{P}$, and $\mathcal{O}$ an acyclic orientation
  of~$\graph{P}$. Then
  \begin{equation}
    \label{eq:walkOrient}
    \#\mathcal{W}\ \leq\ f_2(P)\ \leq\ \Htwosum{\mathcal{O}}
  \end{equation}
  holds. 
  \begin{enumerate}
  \item The first inequality holds with equality if and only
    if~$\mathcal{W}=\twofac{P}$ (i.e., $\mathcal{W}$ ``is'' the set of 
    $2$-faces of~$P$).
  \item The second inequality holds with equality if and only
    if~$\mathcal{O}$ is an AOF-orien\-ta\-tion of~$\graph{P}$.
  \end{enumerate}
\end{theorem}

\begin{proof}
  Since~$\mathcal{O}$ is acyclic, every closed smooth walk in~$\graph{P}$
  must have at least one sink with respect to~$\mathcal{O}$. Thus,
  Lemma~\ref{lem:H2sum} implies
  \begin{equation}
    \label{eq:walkOrientProof1}
    \#\mathcal{W}\leq\Htwosum{\mathcal{O}}\enspace,
  \end{equation}
  yielding
  \begin{equation}
    \label{eq:walkOrientProof2}
    f_2(P)=\#\twofac{P}\leq\Htwosum{\mathcal{O}}\enspace,    
  \end{equation}
  where by Theorem~\ref{thm:2facesAOF} equality holds
  in~\eqref{eq:walkOrientProof2} if and only if~$\mathcal{O}$ is an
  AOF-orientation. Because~$\graph{P}$ has an
  AOF-orientation~$\mathcal{O}_0$ (see Lemma~\ref{lem:manyAOF}),
  inequality~\eqref{eq:walkOrientProof1} gives
  $\#\mathcal{W}\leq\Htwosum{\mathcal{O}_0}=f_2(P)$.  
  
  Hence, it remains to prove that $\#\mathcal{W}=\#\twofac{P}$
  implies~$\mathcal{W}=\twofac{P}$. Suppose, that
  $\#\mathcal{W}=\#\twofac{P}$ holds. It thus suffices to
  show~$\twofac{P}\subseteq\mathcal{W}$ (since we know already
  $\#\twofac{P}\geq\#\mathcal{W}$).  Let~$C\in\twofac{P}$ be any
  closed smooth walk corresponding to a $2$-face of~$P$. By
  Lemma~\ref{lem:AOFterm} there is an AOF-orientation~$\mathcal{O}_C$
  of~$\graph{P}$ such that~$C$ is terminal with respect
  to~$\mathcal{O}_C$. Let~$w_1\in\nodes{P}$ be the unique source
  in~$C$ (with respect to~$\mathcal{O}_C$), and let~$w_0$ and~$w_2$ be
  the two neighbors of~$w_1$ in~$C$. By definition, there is a
  (unique) $W=(w_0,w_1,w_2,\dots,w_{l-1})\in\mathcal{W}$. Because
  of~$\#\mathcal{W}=\#\twofac{P}=\Htwosum{\mathcal{O}_C}$ the closed
  smooth walk~$W$ has a unique sink at some position~$j$ and its
  unique source at position~$1$.  Thus, the two paths
  $(w_1,w_2,\dots,w_j)$ and~$(w_1,w_0,\dots,w_j)$ both are monotone.
  Since~$C$ is terminal, this implies that these two paths are
  contained in~$C$. Therefore we have $C=W\in\mathcal{W}$.
  \begin{figure}[ht]
    \begin{center}
      \epsfig{figure=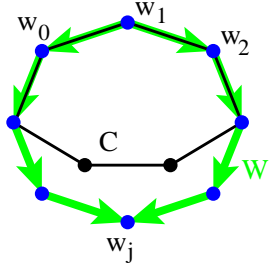,height=3.5cm}
      \caption{Illustration of the proof of Theorem~\ref{thm:walkOrient}. If
        $W\not= C$ then~$C$ cannot be terminal.}
      \label{fig:walkOrient}
    \end{center}
  \end{figure}
\end{proof}

%                                                         
% Good Characterizations                                  
%                                                         
\section{Good Characterizations}
\label{sec:GC}

Theorem~\ref{thm:walkOrient} immediately yields characterizations of
sets of $2$-faces and of AOF-orientations that are similar to
Kalai's characterization of AOF-orientations (see
Lemma~\ref{lem:minHsum}).

\begin{corollary}
  \label{cor:badChar}
  Let~$P$ be a simple polytope. 
  \begin{enumerate}
  \item A facoidal system of walks in~$\graph{P}$ is the
    system~$\twofac{P}$ of $2$-faces of~$P$ if and only if it has
    maximal cardinality among all facoidal systems of walks
    in~$\graph{P}$. 
  \item An acyclic orientation of~$\graph{P}$ is an AOF-orientation if 
    and only if it has minimal $\mathcal{H}_2$-sum among all acyclic
    orientations of~$\graph{P}$. 
  \end{enumerate}
\end{corollary}

Unfortunately, for arbitrary graphs the problem of finding a
$2$-factor with as many cycles as possible is $\NP$-hard. This follows
from the fact that the question whether a graph can be partitioned
into triangles is $\NP$-complete~\cite[Prob.~GT11]{GJ79}.

With respect to algorithmic questions, the following \emph{good
  characterizations} (in the sense of Edmonds~\cite{Edm65b,Edm65a}) of
the set of $2$-faces (and thus, by Proposition~\ref{prop:2facesvfi} of the
vertex-facet incidences) as well as of AOF-orientations may be more
valuable than those in Corollary~\ref{cor:badChar}.

\begin{corollary}
  \label{cor:goodChar}
  Let~$P$ be a simple polytope.
  \begin{enumerate}
  \item Let~$\mathcal{W}$ be a facoidal system of walks
    in~$\graph{P}$. Either there is an acyclic orientation
    of~$\graph{P}$ having a unique sink in every walk
    of~$\mathcal{W}$, or there is a facoidal system of walks
    in~$\graph{P}$ of larger cardinality than~$\#\mathcal{W}$. In the
    first case, $\mathcal{W}=\twofac{P}$ ``is'' the set of $2$-faces
    of~$P$, in the second, it is not.
  \item Let~$\mathcal{O}$ be an acyclic orientation of~$\graph{P}$.
    Either there is a facoidal system~$\mathcal{W}$ of walks
    in~$\graph{P}$ such that~$\mathcal{O}$ has a unique sink in every
    walk in~$\mathcal{W}$, or there is an acyclic orientation
    of~$\graph{P}$ with smaller $\mathcal{H}_2$-sum
    than~$\Htwosum{\mathcal{O}}$. In the first case, $\mathcal{O}$ is
    an AOF-orientation, in the second, it is not.
  \end{enumerate}
\end{corollary}

For graphs~$G$ of simple polytopes let us define \emph{Problem~(A)} as
$$
\max\#\mathcal{W}\quad\text{subject to}\quad\text{$\mathcal{W}$ facoidal 
  system of walks in~$G$}
$$
and \emph{Problem~(B)} as
$$
\min\Htwosum{\mathcal{O}}\quad\text{subject to}\quad\text{$\mathcal{O}$ acyclic
  orientation of~$G$}\enspace.
$$

A third consequence of Theorem~\ref{thm:walkOrient} is the following result.

\begin{corollary}
  \label{cor:minmax}
  The Problems~(A) and~(B) form a pair of strongly dual combinatorial
  optimization problems. The optimal solution of Problem~(A) yields
  the $2$-faces of the respective polytope (and thus its vertex-facet
  incidences, see Proposition~\ref{prop:2facesvfi}).  Every optimal
  solution to Problem~(B) is an AOF-orientation of the graph.
\end{corollary}

Thus, the answer to Perles' original question whether the vertex-facet
incidences of a simple polytope are at all determined by its
 graph, is not only ``yes'' (as proved by Blind and Mani),
or ``yes, and they can be computed'' (as shown by Kalai), but at least
``yes, and they can be computed by solving a combinatorial
optimization problem that has a well-stated strongly dual problem.''

%                                                         
% Remarks                                                 
%                                                         
\section{Remarks}
\label{sec:rem}

Corollary~\ref{cor:minmax} suggests  to design a primal-dual
algorithm for the problem of reconstructing (the vertex-facet
incidences of) a simple polytope from its graph.  Such an
algorithm would start by computing an arbitrary facoidal
system~$\mathcal{W}$ of walks in the given graph (see
Proposition~\ref{prop:facoidal}) and any acyclic orientation~$\mathcal{O}$.
Then it would check for $\#\mathcal{W}=\Htwosum{\mathcal{O}}$. If
equality holds then by Theorem~\ref{thm:walkOrient} one is done.
Otherwise, the algorithm would try to improve either~$\mathcal{W}$
or~$\mathcal{O}$ by exploiting the  reasons for
$\#\mathcal{W}\not=\Htwosum{\mathcal{O}}$. For a concise treatment of
different classical and recent applications of the primal-dual method
in Combinatorial Optimization see~\cite{GW97}. 

Such a (polynomial time) primal-dual algorithm would in particular
yield polynomial time algorithms for the problem to determine an
(arbitrary) AOF-orientation from the graph~$\graph{P}$ of a simple
polytope~$P$ and the set of $2$-faces of~$P$, as well as for the
problem to determine the $2$-faces of~$P$ from~$\graph{P}$ and an
AOF-orientation. As for the first of these two problems it is worth to
mention that no polynomial time method is known that would find any
AOF-orientation even if the input is the entire face lattice of~$P$.
For the second problem no polynomial time algorithm is known as well.

Let (C) be the problem to decide for the graph~$\graph{P}$ of a simple
polytope~$P$ and a set~$\mathcal{C}$ of subsets of nodes
of~$\graph{P}$ if~$\mathcal{C}$ is the set of the subsets of nodes
of~$\graph{P}$ that correspond to the $2$-faces of~$P$. Let (D) be the
problem to decide for the graph~$\graph{P}$ of a simple polytope~$P$
and an orientation~$\mathcal{O}$ of~$\graph{P}$ if~$\mathcal{O}$ is an
AOF-orientation.  The good characterizations in
Corollary~\ref{cor:goodChar} may tempt one to conjecture that these
two problems can be solved in polynomial time.

Unfortunately, from the complexity theoretic point of view,
Corollary~\ref{cor:goodChar} does not provide us with any evidence for
that. In particular, it does \emph{not} imply that problems~(C)
and~(D) are contained in $\NP\cap\coNP$. The reason is that the
problem~(G) to decide for a given graph if there is any simple
polytope~$P$ such that~$G$ is isomorphic to the graph of~$P$ is
neither known to be in~$\NP$ nor in~$\coNP$.
Corollary~\ref{cor:goodChar} only shows that both problems~(C)
and~(D) are in $\NP\cap\coNP$, if one restricts them to any class of
graphs for which problem~(G) is in~$\NP\cap\coNP$.

The problem~(S) to decide for a given lattice~$L$ if there is a simple
polytope~$P$ such that the face lattice~$\lattice{P}$ of~$P$ is
isomorphic to~$L$ is known as the \emph{Steinitz problem for simple
  polytopes}. 

\begin{corollary}
\label{cor:steinitz}
  If problem~(G) is contained in~$\NP$ or~$\coNP$, then problem~(S) is
  contained in~$\NP$ or~$\coNP$, respectively.
\end{corollary}

\begin{proof}
  The face lattice~$\lattice{P}$ of a polytope~$P$ is ranked. The
  graph~$G$ having the rank one elements of~$\lattice{P}$ as its
  nodes, where two rank one elements are adjacent if and only if they
  are below a common rank two element, is isomorphic to the graph
  of~$P$. It can be computed from~$\lattice{P}$ in polynomial time.
  Corollary~\ref{cor:2facesvfi} shows that one can compute the
  vertex-facet incidences (and thus the entire face
  lattice~$\lattice{P}$, see the first paragraph of the introduction)
  of a simple polytope~$P$ in polynomial time (in the size
  of~$\lattice{P}$) from the poset that is induced by the elements of
  rank one (corresponding to the vertices), rank two (corresponding to
  the $1$-faces), and rank three (corresponding to the $2$-faces).
  Together with the first part of Corollary~\ref{cor:goodChar} this
  proves the claim.
\end{proof}

Extending results of Mn\"{e}v~\cite{Mne88} by using techniques
described in~\cite{BS89} one finds (see~\cite[Cor.~9.5.11]{BLSWZ99})
that there is a polynomial (Karp-)reduction of the problem to decide
whether a system of linear inequalities has an integral solution to
problem~(S).  Thus, problem~(S) (and therefore, by
Corollary~\ref{cor:steinitz}, problem~(G)) is not contained in
$\coNP$, unless $\NP=\coNP$.  Furthermore, there are rational simple
polytopes~$P$ with the property that every rational simple
polytope~$Q$ whose graph is isomorphic to the graph of~$P$ has
vertices with super-polynomial coding lengths in the size of the
graphs (this follows from Theorem~B in~\cite{GPS90}).  Thus, it seems
also unlikely that problem~(G) is contained in~$\NP$.

The results presented in Sect.~\ref{sec:WO} hence do neither lead to
efficient algorithms nor to new examples of problems in $\NP\cap\coNP$
not (yet) known to be solvable in polynomial time.  Nevertheless, they
show that the problem to reconstruct a simple polytope from its graph
can be modeled as a combinatorial optimization problem with a
strongly dual problem. We hope that this is an appearance of
Combinatorial Optimization Jack Edmonds is pleased to see in this
volume dedicated to him.

%                                                         
% Acknowledgements                                       
%                                                         
\section{Acknowledgements}

I thank G\"unter M.~Ziegler for valuable comments on an earlier
version of the manuscript.

%                                                         
% References                                              
%                                                         
\bibliographystyle{plain} \bibliography{recsimp}

\begin{thebibliography}{10}

\bibitem{BLSWZ99}
A.~Bj{\"o}rner, M.~Las~Vergnas, B.~Sturmfels, N.~White, and G.~M. Ziegler.
\newblock {\em {Oriented Matroids (2nd ed.)}}, volume~46 of {\em Encyclopedia
  of Mathematics and Its Applications}.
\newblock {Cambridge University Press}, 1999.

\bibitem{BM87}
R.~Blind and P.~Mani-Levitska.
\newblock {Puzzles and polytope isomorphisms.}
\newblock {\em Aequationes Math.}, 34:287--297, 1987.

\bibitem{BS89}
J.~Bokowski and B.~Sturmfels.
\newblock {\em {Computational Synthetic Geometry}}, volume 1355 of {\em Lecture
  Notes in Mathematics}.
\newblock Springer, Heidelberg, 1989.

\bibitem{Dev00}
M.~Develin.
\newblock E-mail conversation, Nov 2000.
\newblock \texttt{develin@bantha.org}.

\bibitem{Edm65b}
J.~Edmonds.
\newblock Maximum matching and a polyhedron with 0,1-vertices.
\newblock {\em J. Res. Natl. Bur. Stand. -- B (Math. and Math. Phys.)},
  69B:125--130, 1965.

\bibitem{Edm65a}
J.~Edmonds.
\newblock Paths, trees, and flowers.
\newblock {\em Can. J. Math.}, 17:449--467, 1965.

\bibitem{GJ79}
M.~R. Garey and D.~S. Johnson.
\newblock {\em {Computers and Intractability. A Guide to the Theory of
  NP-Completeness.}}
\newblock {W. H. Freeman and Company, New York}, 1979.

\bibitem{GW97}
M.~X. Goemans and D.~P. Williamson.
\newblock The primal-dual method for approximation algorithms and its
  application to network design problems.
\newblock In D.~Hochbaum, editor, {\em Approximation Algorithms}, chapter~4.
  PWS Publishing Company, 1997.

\bibitem{GPS90}
J.~E. Goodman, R.~Pollack, and B.~Sturmfels.
\newblock {The intrinsic spread of a configuration in $\R^d$.}
\newblock {\em J. Am. Math. Soc.}, 3(3):639--651, 1990.

\bibitem{HZ01}
C.~Haase and G.~M. Ziegler.
\newblock Examples and counterexamples for {P}erles' conjecture.
\newblock Technical report, TU~Berlin, 2001.
\newblock To appear in: Discrete Comput. Geometry.

\bibitem{HSLdW88}
P.~L. Hammer, B.~Simeone, T.~M. Liebling, and D.~de~Werra.
\newblock From linear separability to unimodality: A hierarchy of
  pseudo-{b}oolean functions.
\newblock {\em SIAM J. Discrete Math.}, 1:174--184, 1988.

\bibitem{JKK01}
M.~Joswig, V.~Kaibel, and F.~K\"orner.
\newblock On the $k$-systems of a simple polytope.
\newblock Technical report, TU~Berlin, 2001.
\newblock arXiv: math.CO/0012204, to appear in: Israel J. Math.

\bibitem{JZ00}
M.~Joswig and G.M. Ziegler.
\newblock {Neighborly cubical polytopes.}
\newblock {\em Discrete Comput. Geometry}, 24:325--344, 2000.

\bibitem{KP01}
V.~Kaibel and M.~Pfetsch.
\newblock Computing the face lattice of a polytope from its vertex-facet
  incidences.
\newblock Technical report, TU~Berlin, 2001.
\newblock arXiv:math.MG/01060043, submitted.

\bibitem{Kal88}
G.~Kalai.
\newblock {A simple way to tell a simple polytope from its graph.}
\newblock {\em J. Comb. Theory, Ser. A}, 49(2):381--383, 1988.

\bibitem{Mne88}
N.~E. Mn{\"e}v.
\newblock {The universality theorems on the classification problem of
  configuration varieties and convex polytopes varieties}.
\newblock In O.Ya. Viro, editor, {\em {Topology and Geometry -- Rohlin
  Seminar}}, volume 1346 of {\em Lecture Notes in Mathematics}, pages 527--543.
  Springer, Heidelberg, 1988.

\bibitem{Tut54}
W.~T. Tutte.
\newblock {A short proof of the factor theorem for finite graphs.}
\newblock {\em Can. J. Math.}, 6:347--352, 1954.

\bibitem{Whi32}
H.~Whitney.
\newblock {Congruent graphs and the connectivity of graphs.}
\newblock {\em Am. J. Math.}, 54:150--168, 1932.

\bibitem{Zie95}
G.~M. Ziegler.
\newblock {\em Lectures on {P}olytopes}.
\newblock Springer-Verlag, New York, 1995.
\newblock Revised edition 1998.

\end{thebibliography}

\end{document}